\documentclass[12pt]{article}
\usepackage[latin1]{inputenc}
\usepackage{amsmath, amsthm, amssymb}
\begin{document}

\bibliographystyle{plain}
\renewcommand{\Im}{\mathop{\rm Im }}
\renewcommand{\Re}{\mathop{\rm Re }}
\newcommand{\ra}{\mathop{\rightarrow }}
\newcommand{\supp}{\mathop{\rm supp}}
\newcommand{\sgn}{\mathop{\rm sgn }}
\newcommand{\card}{\mathop{\rm card }}
\newcommand{\KM}{\mbox{\rm KM}}
\newcommand{\diam}{\mathop{\rm diam}}
\newcommand{\diag}{\mathop{\rm diag}}
\newcommand{\tr}{\mathop{\rm tr}}
\newcommand{\Tr}{\mathop{\rm Tr}}
\newcommand{\dd}{ {\rm d} }
\newcommand{\id}{\mbox{\rm1\hspace{-.2ex}\rule{.1ex}{1.44ex}}
   \hspace{-.82ex}\rule[-.01ex]{1.07ex}{.1ex}\hspace{.2ex}}
\renewcommand{\P}{\mathop{\rm Prob}}
\newcommand{\V}{\mathop{\rm Var}}
\newcommand{\cps}{{\stackrel{{\rm p.s.}}{\longrightarrow}}}
\newcommand{\limm}{\mathop{\rm l.i.m.}}
\newcommand{\cloi}{{\stackrel{{\rm loi}}{\rightarrow}}}
\newcommand{\bra}{\langle\,}
\newcommand{\ket}{\,\rangle}
\newcommand{\obl}{/\!/}
\newcommand{\mapdown}[1]{\vbox{\vskip 4.25pt\hbox{\bigg\downarrow
  \rlap{$\vcenter{\hbox{$#1$}}$}}\vskip 1pt}}
\newcommand{\tab}{&\!\!\!}

\newcommand{\tabb}{&\!\!\!\!\!}
\renewcommand{\d}{\displaystyle}
\newcommand{\epreuve}{\hspace{\fill}$\bigtriangleup$}
\newcommand{\demo}{{\par\noindent{\em D\'emonstration~:~}}}
\newcommand{\solu}{{\par\noindent{\em Solution~:~}}}
\newcommand{\NB}{{\par\noindent{\bf Remarque~:~}}}
\newcommand{\const}{{\rm const}}
\newcommand{\cA}{{\cal A}}
\newcommand{\cB}{{\cal B}}
\newcommand{\cC}{{\cal C}}
\newcommand{\cD}{{\cal D}}
\newcommand{\cE}{{\cal E}}
\newcommand{\cF}{{\cal F}}
\newcommand{\cG}{{\cal G}}
\newcommand{\cH}{{\cal H}}
\newcommand{\cI}{{\cal I}}
\newcommand{\cJ}{{\cal J}}
\newcommand{\cK}{{\cal K}}
\newcommand{\cL}{{\cal L}}
\newcommand{\cM}{{\cal M}}
\newcommand{\cN}{{\cal N}}
\newcommand{\cO}{{\cal O}}
\newcommand{\cP}{{\cal P}}
\newcommand{\cQ}{{\cal Q}}
\newcommand{\cR}{{\cal R}}
\newcommand{\cS}{{\cal S}}
\newcommand{\cT}{{\cal T}}
\newcommand{\cU}{{\cal U}}
\newcommand{\cV}{{\cal V}}
\newcommand{\cW}{{\cal W}}
\newcommand{\cX}{{\cal X}}
\newcommand{\cY}{{\cal Y}}
\newcommand{\cZ}{{\cal Z}}
\newcommand{\bA}{{\bf A}}
\newcommand{\bB}{{\bf B}}
\newcommand{\bC}{{\bf C}}
\newcommand{\bD}{{\bf D}}
\newcommand{\bE}{{\bf E}}
\newcommand{\bF}{{\bf F}}
\newcommand{\bG}{{\bf G}}
\newcommand{\bH}{{\bf H}}
\newcommand{\bI}{{\bf I}}
\newcommand{\bJ}{{\bf J}}
\newcommand{\bK}{{\bf K}}
\newcommand{\bL}{{\bf L}}
\newcommand{\bM}{{\bf M}}
\newcommand{\bN}{{\bf N}}
\newcommand{\bP}{{\bf P}}
\newcommand{\bQ}{{\bf Q}}
\newcommand{\bR}{{\bf R}}
\newcommand{\bS}{{\bf S}}
\newcommand{\bT}{{\bf T}}
\newcommand{\bU}{{\bf U}}
\newcommand{\bV}{{\bf V}}
\newcommand{\bW}{{\bf W}}
\newcommand{\bX}{{\bf X}}
\newcommand{\bY}{{\bf Y}}
\newcommand{\bZ}{{\bf Z}}
\newcommand{\bu}{{\bf u}}
\newcommand{\bv}{{\bf v}}
\newfont{\msbm}{msbm10 scaled\magstep1}
\newfont{\msbms}{msbm7 scaled\magstep1} 
\newcommand{\bbA}{\mbox{$\mbox{\msbm A}$}}
\newcommand{\bbB}{\mbox{$\mbox{\msbm B}$}}
\newcommand{\bbC}{\mbox{$\mbox{\msbm C}$}}
\newcommand{\bbD}{\mbox{$\mbox{\msbm D}$}}
\newcommand{\bbE}{\mbox{$\mbox{\msbm E}$}}
\newcommand{\bbF}{\mbox{$\mbox{\msbm F}$}}
\newcommand{\bbG}{\mbox{$\mbox{\msbm G}$}}
\newcommand{\bbH}{\mbox{$\mbox{\msbm H}$}}
\newcommand{\bbI}{\mbox{$\mbox{\msbm I}$}}
\newcommand{\bbJ}{\mbox{$\mbox{\msbm J}$}}
\newcommand{\bbK}{\mbox{$\mbox{\msbm K}$}}
\newcommand{\bbL}{\mbox{$\mbox{\msbm L}$}}
\newcommand{\bbM}{\mbox{$\mbox{\msbm M}$}}
\newcommand{\bbN}{\mbox{$\mbox{\msbm N}$}}
\newcommand{\bbO}{\mbox{$\mbox{\msbm O}$}}
\newcommand{\bbP}{\mbox{$\mbox{\msbm P}$}}
\newcommand{\bbQ}{\mbox{$\mbox{\msbm Q}$}}
\newcommand{\bbR}{\mbox{$\mbox{\msbm R}$}}
\newcommand{\bbS}{\mbox{$\mbox{\msbm S}$}}
\newcommand{\bbT}{\mbox{$\mbox{\msbm T}$}}
\newcommand{\bbU}{\mbox{$\mbox{\msbm U}$}}
\newcommand{\bbV}{\mbox{$\mbox{\msbm V}$}}
\newcommand{\bbW}{\mbox{$\mbox{\msbm W}$}}
\newcommand{\bbX}{\mbox{$\mbox{\msbm X}$}}
\newcommand{\bbY}{\mbox{$\mbox{\msbm Y}$}}
\newcommand{\bbZ}{\mbox{$\mbox{\msbm Z}$}}
\newcommand{\bbsA}{\mbox{$\mbox{\msbms A}$}}
\newcommand{\bbsB}{\mbox{$\mbox{\msbms B}$}}
\newcommand{\bbsC}{\mbox{$\mbox{\msbms C}$}}
\newcommand{\bbsD}{\mbox{$\mbox{\msbms D}$}}
\newcommand{\bbsE}{\mbox{$\mbox{\msbms E}$}}
\newcommand{\bbsF}{\mbox{$\mbox{\msbms F}$}}
\newcommand{\bbsG}{\mbox{$\mbox{\msbms G}$}}
\newcommand{\bbsH}{\mbox{$\mbox{\msbms H}$}}
\newcommand{\bbsI}{\mbox{$\mbox{\msbms I}$}}
\newcommand{\bbsJ}{\mbox{$\mbox{\msbms J}$}}
\newcommand{\bbsK}{\mbox{$\mbox{\msbms K}$}}
\newcommand{\bbsL}{\mbox{$\mbox{\msbms L}$}}
\newcommand{\bbsM}{\mbox{$\mbox{\msbms M}$}}
\newcommand{\bbsN}{\mbox{$\mbox{\msbms N}$}}
\newcommand{\bbsO}{\mbox{$\mbox{\msbms O}$}}
\newcommand{\bbsP}{\mbox{$\mbox{\msbms P}$}}
\newcommand{\bbsQ}{\mbox{$\mbox{\msbms Q}$}}
\newcommand{\bbsR}{\mbox{$\mbox{\msbms R}$}}
\newcommand{\bbsS}{\mbox{$\mbox{\msbms S}$}}
\newcommand{\bbsT}{\mbox{$\mbox{\msbms T}$}}
\newcommand{\bbsU}{\mbox{$\mbox{\msbms U}$}}
\newcommand{\bbsV}{\mbox{$\mbox{\msbms V}$}}
\newcommand{\bbsW}{\mbox{$\mbox{\msbms W}$}}
\newcommand{\bbsX}{\mbox{$\mbox{\msbms X}$}}
\newcommand{\bbsY}{\mbox{$\mbox{\msbms Y}$}}
\newcommand{\bbsZ}{\mbox{$\mbox{\msbms Z}$}}

%
\def\eurtoday{\number\day \space\ifcase\month\or
 January\or February\or March\or April\or May\or June\or
 July\or August\or September\or October\or November\or December\fi
 \space\number\year}
%
%
\def\aujourdhui{\number\day \space\ifcase\month\or
 janvier\or f{\'e}vrier\or mars\or avril\or mai\or juin\or
 juillet\or ao{\^u}t\or septembre\or octobre\or novembre\or d{\'e}cembre\fi
 \space\number\year}
\textheight=21cm
\textwidth=16cm 
\voffset=-1cm
\newtheorem{theo}{Theorem}[section]
\newtheorem{pr}{Proposition}[section]
\newtheorem{cor}{Corollary}[section]
\newtheorem{lem}{Lemma}[section]
\newtheorem{defn}{Definition}[section]
\hoffset=-1,5cm
\parskip=4mm
 \title{Convergence of dependent walks   in a random 
scenery to  fBm-local time fractional stable motions  }
\author{S.Cohen and C. Dombry   \protect\hspace{1cm}}
\date{}
\maketitle
~
~\\
~\\ 


Serge~{\sc Cohen}, \\
Institut de Mathématique de Toulouse \\
Laboratoire de Statistique et de Probabilit\'{e}s\\
Universit\' e Paul Sabatier\\  
118, Route de  Narbonne, \\
31062 Toulouse  \\
Serge.Cohen@math.univ-toulouse.fr

Clément~{\sc Dombry}, \\
Laboratoire de Mathématiques et Applications (UMR 6086), \\
Université de Poitiers, \\
Téléport 2 - BP 30179, Boulevard Marie et Pierre Curie, 86962 Futuroscope Chasseneuil Cedex\\ clement.dombry@math.univ-poitiers.fr\\

\begin{abstract}
It is classical to approximate the distribution of fractional Brownian
motion by a renormalized sum $ S_n $ of dependent Gaussian random
variables. In this paper we consider  such a walk $ Z_n $ that collects
random rewards $ \xi_j  $ for $ j \in \mathbb Z,$ when the ceiling 
 of the walk $ S_n $ is located at $ j.$ The random reward (or
scenery) $ \xi_j  $ is independent of the walk and with heavy tail.  
We show the convergence of the sum of independent copies of $ Z_n$ suitably renormalized  to a stable 
motion with integral representation, whose kernel is the local time 
of a fractional Brownian motion (fBm). This work extends a previous work where the random walk 
$ S_n$ had independent increments limits. 
\end{abstract}

\section{Introduction} 
\subsection{Motivations}

 Many stochastic processes  have been proposed to model
 communication networks. We can refer to~\cite{WTS95, TWS97, MRRS02} 
for instance, where the limiting processes are either fractional 
Brownian motion or L\'evy $\beta$-stable process. More recently, in
 \cite{CS05} a  process named $H$-fBm local time fractional stable 
motion was constructed.  When $H=\frac{1}{2}$,   the so called Random Reward 
 Schema, was  also proposed, it is  a  discrete schema, which could be thought  of as a toy model for
 INTERNET traffic, and which is converging to  this process. The aim of this paper is 
to extend these results to the case $H\neq \frac{1}{2}$. In the proof
of the convergence in~\cite{CS05} a strong approximation of the 
local time of standard Brownian motion was used. As far as we know no 
strong approximation of the local time of fractional Brownian motion 
is available and it was one the problems to overcome. In~\cite{DG07}  
discrete approximations of local time fractional stable 
motion have been obtained by Dombry and Guillotin-Plantard
\cite{DG07} where the fBm local time is replaced by the local time of
an $\alpha$-stable Levy motion. But they did not use strong
approximation of the local time. In this paper we use  convergence of the local time of a 
classical walk with dependent increments to the  local time of
fractional Brownian motion and the technique in~\cite{DG07} to get our result. 
Please note that other approximations of fBm local time fractional stable
motion have been considered in~\cite{Marouby08}, but they are not
related to walks in random sceneries. 

\subsection{Model and results}
Let $\xi=(\xi_x)_{x\in\bbsZ}$ denote a sequence of independent, identically distributed, symmetric real-valued random variables. The sequence $\xi$ is called a {\it random scenery}. Suppose that it belongs to the normal domain of attraction  of a stable symmetric distribution $Z_{\beta}$ of index $\beta\in (0,2]$. This means that the following weak convergence holds:
\begin{equation}\label{eq1.01}
n^{-\frac{1}{\beta}}\sum_{x=0}^n \xi_x \mathop{\Longrightarrow}_{n\rightarrow\infty}^{\cL} Z_\beta, 
\end{equation}
where $Z_\beta$ is the symmetric stable law with characteristic function $\bar\lambda$ given by
\begin{equation}\label{eq1.02}
\bar\lambda(u)=\bbE\exp(iuZ_\beta)=\exp\left(-\sigma^\beta|u|^\beta \right),\ \ u\in\bbR 
\end{equation}
for some constants $\sigma>0$. 

Let $S=(S_k)_{k\in\bbsN}$ be a {\it random walk}  on $\bbZ$ 
independent of the random scenery $\xi$. We suppose that 
\begin{equation}\label{eq1.04}
\left\{\begin{array}{l}S_0=0, \\ 
S_n=\sum_{k=1}^n X_k\ ,\ n\geq 1,\end{array}\right.
\end{equation}
where $X_i,i\geq 1$ is a stationary Gaussian sequence with mean $0$ and correlations $r(i-j)=\bbE[X_iX_j]$ satisfying
\begin{equation}\label{eq1.05}
\sum_{i=1}^n\sum_{j=1}^n r(i-j)\sim n^{2H},
\end{equation}
as $n\rightarrow\infty$, with $0<H<1$. \\

We define the {\it random walk in random scenery} as the process $(Z_n)_{n\geq 0}$ given by
\begin{equation}\label{eq1.06}
Z_n=\sum_{k=0}^n \xi_{([S_k])},
\end{equation}
where $[S_k]$ is the ceiling of $S_k$.
Stated simply, a random walk in random scenery is a cumulative sum process whose summands are drawn from the scenery; the order in which the summands are drawn is determined by the path of the random walk.
We extend this definition to non-integer time $s\geq 0$ by the linear interpolation
\begin{equation}\label{eq1.07}
Z_s=Z_{[s]}+(s-[s])(Z_{[s]+1}-Z_{[s]}).
\end{equation}

We now describe the limit theorem for the random walk in random
scenery established by Wang \cite{W03} 
(in the case $\beta=2$).\\ 
Cumulative sums of the scenery converge in $\cD(\bbR)$, the space of càd-làg functions:
$$\left(n^{-\frac{1}{\beta}} \sum_{k=0}^{[nx]}\xi_k\right)_{x\in \bbsR} \mathop{\Longrightarrow}_{n\rightarrow\infty}^{\cL} \left(W(x)\right)_{x\in \bbsR},$$ 
where $W$ is a bilateral $\beta$-stable Lévy process such that $W(0)=0$, and $W(1)$ and $W(-1)$ are distributed according to  $Z_\beta$.\\
The covariance structure of the sequence $X_i$ given by equation (\ref{eq1.05}) implies that $S_n,n\geq 0$ belongs to the domain of attraction of the fractional Brownian motion of Hurst index $H$, i.e. the following convergence hold in $\cD([0,\infty))$ (cf.~\cite{T79}.) 
\begin{equation}\label{1.08}
\frac{1}{n^{H}} (S_{[nt]})_{0\leq t\leq 1} \mathop{\Longrightarrow}_{n\rightarrow\infty}^{\cL} B_H(t), 
\end{equation}
To describe the limit process known as {\it fractional Brownian motion
  in stable scenery} we suppose that $B_H$ and $W$ are two independent
processes defined on the same probability space and distributed as
above. Let $L_t(x)$ the jointly continuous version of the local time
of the process $B_H$ (cf.~\cite{B73}).

In the case $\beta=2$ corresponding to the case of a Gaussian scenery, Wang proves the following weak convergence in the space of continuous function $\cC([0,\infty))$  
\begin{equation}\label{eq1.09}
\left(n^{-\delta} Z_{nt}\right)_{t\geq 0} \mathop{\Longrightarrow}_{n\rightarrow\infty}^{\cL}  \left(\Delta(t)\right)_{t\geq 0}
\end{equation} 
where $\delta=1-H+H\beta^{-1}$ and $\Delta$ is the process defined by
$$\Delta(t)=\int_{-\infty}^{+\infty} L_t(x){\rm d}W(x).$$
The limit process $\Delta$ is a  continuous $\delta$-self-similar stationary increments process.


Our results state the convergence of the so called Random Reward Schema to the fBm local time stable fractionnal motion.
We begin with a continuous version and consider $\Delta^{(i)}, i\geq 1$ independent copies of the process $\Delta$. 
\begin{theo}\label{theo1} 
The following weak convergence holds in $\cC([0,\infty))$:
\begin{equation}\label{eq1.10}
\left(n^{-\frac{1}{\beta}} \sum_{i=1}^n \Delta^{(i)}(t)\right)_{t\geq 0} \mathop{\Longrightarrow}_{n\rightarrow\infty}^{\cL}  \left(\Gamma(t)\right)_{t\geq 0},
\end{equation}
where $\Gamma$ is a H-fBm local time stable fractionnal motion.
\end{theo}
\noindent
Replacing the stable process in random scenery by a random walk in random scenery, we obtain the random rewards schema which yields a discrete approximation of the process $\Gamma$.
Let $\xi^{(i)}=(\xi^{(i)}_x)_{x\in \bbsZ},\ i\geq 1$ be independent copies of $\xi$. Let $S^{(i)}=(S^{(i)}_n)_{n\in\bbsN}$ be independent copies of $S$ and also independent of the $\xi^{(i)}, i\geq 1$.
Denote by $D_n^{(i)}$ the $i$-th random walk in random scenery defined by 
\begin{equation}\label{eq1.11}
D_n^{(i)}(t)=n^{-\delta} Z_{nt}^{(i)} 
\end{equation}
where the definition of $Z_n^{(i)}$ is given by equations (\ref{eq1.06}) and (\ref{eq1.07}) with $\xi$ and $S$ replaced by the $i$-th random scenery $\xi^{(i)}$ and the $i$-th random walk $S^{(i)}$ respectively. 
\begin{theo}\label{theo2}
Let $c_n$ be a sequence of integers such that $\lim c_n=+\infty$. Then, the following weak convergence holds in $\cC([0,\infty))$:
\begin{equation}\label{eq1.12}
\left(c_n^{-\frac{1}{\beta}} \sum_{i=1}^{c_n} D_n^{(i)}(t)\right)_{t\geq 0} \mathop{\Longrightarrow}_{n\rightarrow\infty}^{\cL}  \left(\Gamma(t)\right)_{t\geq 0}.
\end{equation}
The limit process $\Gamma$ is the same as in Theorem \ref{theo1}.
\end{theo}

\section{Sums of stable processes in random scenery \\ Proof of Theorem \ref{theo1}}

For $n\geq 1$, let $\Gamma_n$ the continuous process defined by 
$$\Gamma_n(t)=n^{-\frac{1}{\beta}}\sum_{i=1}^n \Delta^{(i)}(t)\ ,\ t\geq 0.$$
Theorem \ref{theo1} claims that the sequence $\Gamma_n$ converges weakly in $\cC([0,\infty))$. 
We prove this fact by proving the convergence of the finite dimensional distributions and the tightness of the sequence.
Theorem \ref{theo1} is thus a consequence of Propositions \ref{pr1.1} and \ref{pr1.2} below. \\

We first need a Lemma giving the characteristic function of the finite dimensional distribution of $\Delta$:
\begin{lem}\label{lem1.1}
For any $(\theta_1,\cdots,\theta_k)\in \bbR^k$ and $(t_1,\cdots,t_k)\in [0,+\infty)^k$
$$\bbE\left[\exp\left( i\sum_{j=1}^k \theta_j \Delta(t_j) \right)\right]=\bbE\left[\exp(-\sigma^{\beta}X)\right]$$
with
\begin{equation}
  \label{eq:X}
 X= \int_{\bbsR} |\sum_{j=1}^k \theta_jL_{t_j}(x)|^\beta {\rm d}x .
\end{equation}
\end{lem}
{\it Proof :} This is the analogous of Lemma 5 in Kesten and Spitzer
giving the characteristic function of the finite dimensional
distribution of the stable Levy-process in stable scenery, for the
fractional Brownian motion in stable scenery. 
The demonstration is the same replacing the local time of a stable Lévy process by the local time of the fractional Brownian motion.  $\Box$\\

\begin{pr}\label{pr1.1}
The finite dimensional distributions of $(\Gamma_n(t))_{t\geq 0}$ converge weakly as $n\rightarrow\infty$ to those of $(\Gamma(t))_{t\geq 0}$. 
\end{pr} 
{\it Proof :} \ Let $(\theta_1,\cdots,\theta_k)\in \bbR^k$ and $(t_1,\cdots,t_k)\in [0,+\infty)^k$. 
We compute the characteristic functions 
\begin{eqnarray}
\bbE\left[\exp\left( i\sum_{j=1}^k \theta_j\Gamma_n(t_j)\right)\right]&=&\bbE\left[\exp\left( in^{-\frac{1}{\beta}}\sum_{j=1}^k \theta_j \Delta(t_j) \right)\right]^n \nonumber\\
&=&\bbE\left[\exp(-n^{-1}X)\right]^n \label{eq2.1}
\end{eqnarray}
We prove that the following asymptotic holds:
\begin{equation}\label{eq2.2}
\bbE\left[\exp(-n^{-1}\sigma^{\beta}X)\right]=1-n^{-1}\sigma^{\beta}\bbE(X)+o(n^{-1}).
\end{equation}
Note that the integrability of the random variable $X$ follows from the inequality
$$|X|\leq \left(\sum_{j=1}^k |\theta_j|\right)^\beta \int_{\bbsR} L_t(x)^{\beta}{\rm d}x$$
where $t=\max\{t_j\ ,\ 1\leq j\leq k\}$, and the fact that $\bbE\int_{\bbsR} L_t(x)^{\beta}{\rm d}x<\infty $ which is proved in \cite{CS05} Theorem 3.1. \\
We now prove equation (\ref{eq2.2}). To this aim, observe that
$$n\left(\bbE\left[\exp(-n^{-1}\sigma^\beta X)\right]^n-1\right)=\bbE(f_n(X))\mathop{\longrightarrow}_{n\rightarrow\infty} -\sigma^\beta \bbE(X)$$
where $f_n$ is defined on $\bbC$ by $f_n(x)=n(\exp(-n^{-1}\sigma^\beta x)-1)$. The convergence follows from the dominated convergence Theorem because $f_n(X)$ converges almost surely to $-\sigma^\beta X$ and $|f_n(X)|$ is almost surely bounded from above by $\sigma^\beta |X|$ which is integrable.
Finally, equations (\ref{eq2.1}) and (\ref{eq2.2}) together yield 
$$\bbE\left[\exp\left( i\sum_{j=1}^k \theta_j\Gamma_n(t_j)\right)\right]\mathop{\longrightarrow}_{n\rightarrow\infty}\exp\left(-\sigma^\beta\bbE(X)\right).$$
This proves Proposition \ref{pr1.1}.
$\Box$\\

\begin{pr}\label{pr1.2}
The sequence of process $\Gamma_n$ is tight in $\cC([0,\infty))$.
\end{pr} 
{\it Proof : } We follow the proof of Proposition 2.2 in \cite{DG07} and give only the main lines of the proof, the details are to be found in\cite{DG07}. The difference is that the $\alpha$-stable Levy motion $Y_t$ is replaced by the fBm $B_H(t)$ of index $H.$ Hence the local time process of $Y$ is replaced by the local time of $B_H$ and denoted in both context by $L_t(x)$. Furthermore, the self-similarity index of $Y$ is equal to $1/\alpha$ and has to be replaced by $H$. \\
The case $\beta=2$ is straightforward and relies on Itô's isometry: the process $\Gamma_n$ is square integrable and for all $0\leq t_1<t_2$
\begin{eqnarray*}
\bbE\left[|\Gamma_n(t_2)-\Gamma_n(t_1)|^2\right]&=& \bbE\left[|n^{-\frac{1}{2}} \sum_{i=1}^n \Delta^{(i)}(t_2)-\Delta^{(i)}(t_1) |^2\right]\\ 
  &=&  \sigma^2(t_2-t_1)^{2-H}\bbE\left[\int_{\bbsR} L_{1}(x)^2\dd x\right]. 
\end{eqnarray*}
Using Kolmogorov criterion, we deduce that the sequence $\Gamma_n$ is tight. 

In the case $0<\beta<2$, the process $\Gamma_n$ has infinite variance and we use the truncation method. 
Introduce the Lévy-Itô decomposition of $W$: 
\begin{equation}\label{eq2.3}
W_x=bx+\int_0^x\int_{|u|\leq 1} u(\mu-\bar\mu)(\dd u,\dd s)+\int_0^x\int_{|u|> 1} u\mu (\dd u,\dd s)
\end{equation}
where $b$ is the drift and $\mu$ is a Poisson random measure on $\bbR\times\bbR$ with intensity $\bar\mu(\dd u,\dd x)=\lambda(\dd u)\otimes \dd x$,  and $\lambda$ is the stable Lévy measure on $\bbR$:
$$\lambda(\dd u)=\left(c_-1_{\{u<0\}}+c_+1_{\{u>0\}}\right) \frac{\dd u}{|u|^{\beta +1}},\ \ c_-,c_+\geq 0, c_-+c_+>0.$$
For some truncation level $R>1$, let $W^{(R^-)}$ and $W^{(R^+)}$ be the independent Lévy processes defined by 
$$W_x^{(R^-)}=\int_0^x\int_{|u|\leq R} u(\mu-\bar\mu)(\dd u,\dd s),\ \ \ W_x^{(R^+)}=\int_0^x\int_{|u|> R} u\mu (\dd u,\dd s).$$
The Lévy-Itô decomposition (\ref{eq2.3}) rewrites as 
$$W_x=b_R x+W_x^{(R^-)}+W_x^{(R^+)}$$
where $b_R=b+\int_{1<|y|\leq R} u\lambda(\dd u)$ is a drift depending on $R$.
This decomposition of the stable scenery yields the following decomposition of the stable process in random scenery:
$$\Delta(t)=b_Rt + \Delta^{(R^-)}(t) + \Delta^{(R^+)}(t),$$
with 
$$\Delta^{(R^-)}(t)=\int_{\bbsR} L_t(x)\dd W_x^{(R^-)}, \ \ \ \Delta^{(R^+)}(t)=\int_{\bbsR} L_t(x)\dd W_x^{(R^+)} .$$
Considering such decomposition for i.i.d. copies of $\Delta$, we have with obvious notations, the following decomposition of $\Gamma_n$: 
\begin{equation}\label{eq2.9}
\Gamma_n(t)=n^{1-\frac{1}{\beta}} b_{R_n} t +\Gamma_n^{(R^-)}(t)+\Gamma_n{(R^+)}(t),
\end{equation}
with 
$$\Gamma_n^{(R^-)}(t)=n^{-\frac{1}{\beta}}\sum_{i=1}^n \Delta^{(i,R_n^-)}(t), \ \ \Gamma_n^{(R^+)}(t)=n^{-\frac{1}{\beta}}\sum_{i=1}^n \Delta^{(i,R_n^+)}(t)$$
with truncation level $R_n=Rn^{\frac{1}{\beta}}$.
The sequence $n^{1-\frac{1}{\beta}} b_{R_n}$ is known to be bounded (assertion A1 in \cite{DG07}).\\
Similarly to equation (23) in \cite{DG07}, the process $\Gamma_n^{(R^-)}(t)$  is square integrable and for any $0\leq t_1<t_2$,
\begin{eqnarray}
\bbE\left[(\Gamma_n^{(R^-)}(t_2)- \Gamma_n^{(R^-)}(t_1))^2\right] &=&  \frac{c_-+c_+}{2-\beta} R^{2-\beta} (t_2-t_1)^{2-H}\bbE\left[\int_{\bbsR} L_{1}(x)^2\dd x\right]\nonumber.
\end{eqnarray}
Using Kolmogorov criterion, this estimate implies that the sequence of process $\Gamma_n^{(R^-)}$ is tight.\\
On the other hand, similarly to assertion (A3) in \cite{DG07}, the probability that $\Gamma_n^{(R^+)}\equiv 0 $ on $[0,T]$ satisfies
\begin{eqnarray*}
\bbP\left( \Gamma_n^{(R^+)}\equiv 0 {\rm\ on\ }[0,T] \right)&\geq& \left[\bbP\left( \Delta^{(R_n^+)}\equiv 0 {\rm\ on\ }[0,T] \right)\right]^n \\
 &\geq& \left[  1-2\frac{c^++c^-}{\beta}R_n^{-\beta} \bbE\left(\sup_{0\leq t\leq T} |B_H(t)| \right)\right]^n.
\end{eqnarray*}
and hence 
$$\lim_{R\rightarrow\infty}\limsup_{n\rightarrow\infty} \bbP\left(\Gamma_n^{(R^+)}\equiv 0 {\rm\ on\ }[0,T]\right)=1.$$
These facts imply the tightness of the sequence $\Gamma_n$.
$\Box$\\

\section{Fractional Random Reward Schema \\ Proof of Theorem \ref{theo2}}
We define the process $G_n$ by 
\begin{equation}\label{eq3.14}
G_n(t)= c_n^{-\frac{1}{\beta}} \sum_{i=1}^{c_n} D_n^{(i)}(t),\ \ t\geq 0,
\end{equation}
where $D_n^{(i)}$ is the i-th random walk in random scenery properly rescaled and defined by (\ref{eq1.11}). 
Theorem \ref{theo2} states that $G_n$ converges weakly to $\Gamma$ in $\cC([0,\infty))$. The key tool in the proof  is the local time of the strongly correlated random walk $(S_k)_{k\geq 0}$ (we omit the superscript $(i)$).\\

Let $x\in\bbZ$ and $n\geq 1$. The local time $N_n(x)$ of the random walk $(S_k)_{k\geq 0}$ at point $x$ up to time $n$ is defined by
$$N_n(x)=\sum_{k=0}^n 1_{\{[S_k]=x\}}.$$
It represents the amount of time the walk spends in the interval $[x,x+1[$ up to time $n$. 
We extend this definition to non-integer time $s\geq 0$ by linear interpolation:
$$N_s(x)=N_{[s]}(x)+(s-[s])(N_{[s]+1}(x)-N_{[s]}(x)).$$
The random walk in random scenery  writes for all $s\geq 0$  
\begin{equation}\label{eq3.1}
D_n(t)=n^{-\delta}\sum_{x\in \bbsZ} N_{nt}(x)\xi_x 
\end{equation}
where the collection of random variables $\{N_s(x), x\in\bbZ\}$ and $\{\xi_x, x\in\bbZ\}$ are independent.\\

We collect in the next subsection different results about the local
times of the strongly correlated random walks that will be of great
use in the sequel. 
Although the results are analogous to the ones in \cite{DG07} for independent increments random walks, some difficulties arise from the strong correlations of the increments. However in \cite{W03}, Wang shows how to use the Gaussian structure to get some estimates on the local times of the strongly correlated random walk.

Then the proof of Theorem \ref{theo2} is quite analogous to the proof of Theorem 2 in \cite{DG07}. Proposition \ref{pr2.1} states the convergence of the finite dimensional distribution. The tightness of the sequence is stated in Proposition \ref{pr2.2}. We give the main lines of the proof and omit some details that are to be found in \cite{DG07}.

\subsection{Some results about local times}

\subsubsection{Maximum local time, self intersection local time and range.}
The maximum local time $L_n$ of the random walk up to time $n$ is defined by
$$L_n=\sup_{x\in\bbsZ} N_n(x).$$
The number of self-intersections $V_n$ of the random walk up to time $n$ is defined by
$$V_n=\sum_{0\leq i,j\leq n} 1_{\{[S_i]=[S_j]\}}=\sum_{x\in\bbsZ}N_n(x)^2.$$
The range $R_n$ of the random walk up to time $n$ is defined by
$$ R_n=\sum_{x\in\bbsZ}1_{\{N_n(x)\neq 0\}}.$$
These definitions extend obviously to non-integer time $s\geq 0$.

Our results rely on different estimations of these quantities that we gather in the following Lemma:
\begin{lem}\label{lem3.1}\ \\
\vspace{-0.3cm}
\begin{itemize}
\item The following convergence in probability holds 
\begin{equation}\label{eq3.4}
n^{-\delta}L_n \mathop{\longrightarrow}_{n\rightarrow\infty}^{{\rm P}} 0.
\end{equation} 
\item  For any $p\in [1,+\infty)$, there exists some constant $C$ such that for all $n\geq 1$, 
\begin{equation}\label{eq3.5}
\bbE\left(V_n^p\right)\leq C n^{p(2-H)}.
\end{equation}
\item For any $p\in [1,+\infty)$, there exists some constant $C$ such that for all $n\geq 1$
\begin{equation}\label{eq3.6}
\bbE\left( R_n^p\right)\leq Cn^{pH}.
\end{equation}
\end{itemize}
\end{lem}
\noindent

{\it Proof:}\ \\
$\bullet$ We follow the lines of Lemma 4 in Kesten and Spitzer. Let $\varepsilon>0$, we have,
\begin{eqnarray*}
\bbP(n^{-\delta}L_n>\varepsilon )&\leq & \bbP(N_n(x)>0 {\rm\ for\ some\ } |x|> An^H )+\sum_{|x|\leq An^H} \bbP(N_n(x)>n^\delta\varepsilon)\\
&\leq & \bbP(\sup_{0\leq k\leq n} n^{-H}|[S_k]|>A)+\sum_{|x|\leq An^H} \bbE(N_n(x)^p)n^{-p\delta}\varepsilon^{-p}
\end{eqnarray*}
We now use the following estimation from \cite{W03} lemma 4.4: there exists some $C>0$ such that
$$\bbE(N_n(x)^p)\leq Cn^{p(1-H)},$$
and hence we have for all $A>0$ and $\varepsilon>0$
$$\sum_{|x|\leq An^H} \bbE(N_n(x)^p)n^{-p\delta}\varepsilon^{-p}\leq 2An^HCn^{p(1-H)}n^{-p\delta}\varepsilon^{-p},$$
and this quantity goes to $0$ as $n\rightarrow\infty$ if we choose $p$ large enough such that $H+p(1-H)-p\delta=H(1-p/\beta)<0$.
At last, the term $$\bbP(\sup_{0\leq k\leq n} n^{-H}|[S_k]|>A)$$ converges to $\bbP(\sup_{0\leq t\leq 1} |B_H(t)|>A)$ as $n\rightarrow\infty$, and this last term goes to zero as $A\rightarrow\infty$.\\

$\bullet$  First notice that we can suppose without restriction that $p\geq 1$ is an integer, because the bound for $p'\geq 1$ is a consequence of the case  $p\geq p'$.\\
The number of self-intersections up to time $n$ is bounded from above
by
$$V_n\leq \sum_{0\leq i\leq j\leq n}2 1_{\{[S_i]=[S_j]\}}.$$
Using Minkowski inequality,
\begin{equation}\label{eq3.7}
|\!|V_n|\!|_p \leq  2\sum_{i=0}^n |\!| \sum_{j=i}^n 1_{\{[S_i]=[S_j]\}}|\!|_p,
\end{equation}
where $|\!|X|\!|_p=\bbE(|X|^p)^{1/p}$.
For fixed $i$, the stationarity of the random walk's increments implies that the distribution of $\sum_{j=i}^n 1_{\{[S_i]=[S_j]\}}$ and 
$\sum_{j=0}^{n-i} 1_{\{[S_i]=0\}}=N_{n-i}(0)$ are equal. 
Since $N_{n-i}(0)\leq N_n(0)$, equation (\ref{eq3.7}) yields
\begin{equation}\label{eq3.8}
\bbE(V_n^p)\leq  2^pn^p \bbE\left(N_n(0)^p\right).
\end{equation}
We now refer to lemma 4.4 in \cite{W03} which states that there is some $C>0$ such that
\begin{equation}\label{eq3.9}
\bbE\left(N_n(0)^p\right)\leq C n^{p(1-H)}
\end{equation}
Equations (\ref{eq3.8}) and (\ref{eq3.9}) together yield equation (\ref{eq3.5}).\\
$\bullet$  We only have to notice that $R_n\leq 1+2\sup_{0\leq k\leq n} |S_k|$ and hence it is enough to prove that $\sup_{0\leq k\leq n} n^{-H}|S_k|$ is bounded in $L^p$ for all $p\geq 1$. Let $S^n(t)_{t\in [0 1]}$ be the continuous process defined by 
$$S^n(t)=n^{-H}S_{[nt]} +(nt-[nt])n^{-H}(S_{[nt]+1}-S{[nt]}).$$
By equation (\ref{eq1.05}), the sequence of process $S^n$ converges weakly to $B_H$ in $\cC([0,1])$ furnished with the uniform norm $|\!|.|\!|_\infty$. Furthermore, 
$$\sup_{0\leq k\leq n} n^{-H}|S_k|=|\!|S^n|\!|_\infty.$$
 Hence we need to show that $|\!|S^n|\!|_\infty$ is bounded in $L^p$ for all $p\geq 1$. Using a concentration result (see e.g. \cite{LT91} p. 60), a sequence of Gaussian random variables which is bounded in probability is bounded in all $L^p$ spaces. Since the sequence $S_n$ converges in distribution to $B_H$, it is bounded in probability, and hence bounded in all $L^p$ spaces.

\subsubsection{Convergence of functional of local times}
The following lemma is an analogous of Lemma 6 in \cite{KS79} when
 the random walk in the domain of attraction of a stable Levy motion
 is replaced by a random walk in the domain of attraction of a
 fractional Brownian motion Note that the case of 
walks in a $L^2$ scenery was considered by Wang. We here generalize
Wang's result (Proposition 3.2 in \cite{W03}) to the case of an heavy
tailed scenery.

\begin{lem}\label{KS}
For all $(\theta_1,\cdots,\theta_k)\in\bbR^k$, $(t_1,\cdots,t_k)\in [0,+\infty)^k$, $\sigma>0$, $\beta\in(0,2]$, the distribution of
$$X_n= n^{-\delta\beta}\sum_{x\in\bbsZ}\big|\sum_{j=1}^k \theta_j N_{nt_j}(x)\big|^\beta  $$
converges weakly as $n\rightarrow\infty$ to $X$ defined by equation 
                                (\ref{eq:X}).\\
Furthermore, $X_n$ is bounded in $L^p$ for all $p\geq 1$.
\end{lem}
{\it Proof:}\\
Following Kesten and Spitzer \cite{KS79} and Wang \cite{W03}, we introduce for small $\tau>0$ and large $N$,
$$U(\tau,N,n)=n^{-\beta\delta} \sum_{|x|\leq N\tau n^H }\big|\sum_{j=1}^k \theta_j N_{nt_j}(x)\big|^\beta $$
$$d(l,n)=n^{-1}\sum_{j=1}^k \theta_j \sum_{l\tau n^H\leq y < (l+1)\tau n^H}   N_{nt_j}(y)$$
$$V(\tau,N,n)=\tau^{1-\beta}\sum_{|l|\leq N} |d(l,n)|^\beta .$$
Then, 
\begin{eqnarray}
& &X_n-U(\tau,N,n)-V(\tau,N,n)\nonumber \\
&=&\sum_{|l|\leq N} \sum_{l\tau n^H\leq x < (l+1)\tau n^H}n^{-\delta\beta}\left\{ \big|\sum_{j=1}^k \theta_j N_{nt_j}(x)\big|^\beta -n^\beta[\tau n^H]^{-\beta}|d(l,n)|^\beta\right\}\nonumber \\
& &+ \sum_{|l|\leq M}(n^{\beta-\delta\beta}[\tau n^{H}]^{1-\beta}-\tau^{1-\beta}) |d(l,n)|^\beta. \label{eqq}
\end{eqnarray}
By Lemma 3.1 in \cite{W03} $d(l,n)$ converges in distribution to 
$$\sum_{j=1}^k \theta_j \int_{l \tau}^{(l+1)\tau} L_{t_j}(x)\dd x.$$
Since 
$$ n^{\beta-\delta\beta}[\tau n^{H}]^{1-\beta}-\tau^{1-\beta}\rightarrow 0,$$ 
the second sum over $l$ in the right hand side of (\ref{eqq}) tends to zero in probability as $n\rightarrow \infty$.
We now show that the first sum over $l$ in the right hand side of (\ref{eqq}) is small in probability when $\tau$ is small. We use the following inequality, valid for any $a\geq 0$, $b\geq 0$
$$ |a^\beta-b^\beta|\leq \left\{ \begin{array}{lll} |a-b|^\beta & {\rm if} & \beta\leq 1 \\ \beta|a-b|(a^{\beta-1}+b^{\beta-1}) & {\rm if} & \beta<1,\end{array}\right.$$
to estimate the sum over $x$. In the case $\beta\leq1$,  
\begin{eqnarray*}
 & &\bbE \left| \sum_{|l|\leq N} \sum_{l\tau n^H\leq x< (l+1)\tau n^H}n^{-\delta\beta}\left\{ \big|\sum_{j=1}^k \theta_j N_{nt_j}(x)\big|^\beta -n^\beta[\tau n^H]^{-\beta}|d(l,n)|^\beta\right\} \right|\\
&\leq& n^{-\delta\beta} \sum_{|l|\leq N}\sum_{l\tau n^H\leq x< (l+1)\tau n^H} \left[ \bbE \left|  \sum_{j=1}^k \theta_j N_{nt_j}(x) -n[\tau n^H]^{-1}d(l,n) \right|^2 \right]^{\beta/2}
\end{eqnarray*}
and 
\begin{eqnarray*}
 & &\bbE \left|  \sum_{j=1}^k \theta_j N_{nt_j}(x) -n[\tau n^H]^{-1}d(l,n) \right|^2\\
&\leq& [\tau n^H]^{-1}\sum_{i=1}^k \theta_i^2 \sum_{j=1}^k \sum_{l\tau n^H\leq y<(l+1)\tau n^H } \bbE|N_{nt_j}(x)-N_{nt_j}(y)|^2.
\end{eqnarray*}
Now from Lemma 4.6 in \cite{W03}, there exist $C$ and $r>0$ such that for large $n$, large $N$  and small $\tau$ (with $A=N\tau$ large enough), and any $|x|\leq A$ and $y$ such that $|x-y|\leq \tau n^H$, the following holds : 
$$\bbE|N_{nt_j}(x)-N_{nt_j}(y)|^2 \leq C A \tau^r n^{2-2H} .$$
Combining these estimates,  
 \begin{eqnarray*}
 & &\bbE \left| \sum_{|l|\leq N} \sum_{l\tau n^H\leq x< (l+1)\tau n^H}n^{-\delta\beta}\left\{ \big|\sum_{j=1}^k \theta_j N_{nt_j}(x)\big|^\beta -n^\beta[\tau n^H]^{-\beta}|d(l,n)|^\beta\right\} \right|\\
&\leq& n^{-\delta\beta}(2N+1)[n^H\tau]\left([n^H\tau]^{-1}[n^H\tau]CA\tau^rn^{2-2H}\right)^{\beta/2} \\
&\leq& C N\tau^{1+r\beta/2}
\end{eqnarray*}
This completes the estimate of (\ref{eqq}) in the case $\beta\leq 1$.

In the case $\beta>1$, we have
\begin{eqnarray*}
 & &\bbE \left| \sum_{|l|\leq N} \sum_{l\tau n^H\leq x< (l+1)\tau n^H}n^{-\delta\beta}\left\{ \big|\sum_{j=1}^k \theta_j N_{nt_j}(x)\big|^\beta -n^\beta[\tau n^H]^{-\beta}|d(l,n)|^\beta\right\} \right|\\
&\leq& n^{-\delta\beta} \sum_{|l|\leq N}\sum_{l\tau n^H\leq x< (l+1)\tau n^H} \beta \bbE \left[\left|  \sum_{j=1}^k \theta_j N_{nt_j}(x) -n[\tau n^H]^{-1}d(l,n) \right|\right.\\
& &\left.\times \left( |\sum_{j=1}^k \theta_j N_{nt_j}(x)|^{\beta-1}+|n[\tau n^H]^{-1}d(l,n)|^{\beta-1}\right)\right]\\
& \leq&\beta n^{-\delta\beta}  \left[\sum_{|l|\leq N}\sum_{l\tau n^H\leq x< (l+1)\tau n^H}\bbE \left|  \sum_{j=1}^k \theta_j N_{nt_j}(x) -n[\tau n^H]^{-1}d(l,n) \right|^2\right]^{1/2}\\
& &\times \left[\sum_{|l|\leq N}\sum_{l\tau n^H\leq x< (l+1)\tau n^H}\bbE \left( |\sum_{j=1}^k \theta_j N_{nt_j}(x)|^{\beta-1}+|n[\tau n^H]^{-1}d(l,n)|^{\beta-1}\right)^2\right]^{1/2}.
\end{eqnarray*}
With the same techniques as above, the first factor is shown to be bounded from above by
$$\beta n^{-\delta\beta}[(2N+1)\tau n^HCA\tau^rn^{2-2H}]^{1/2}$$
In order to 
upper-bound  the second factor, introduce $T=\sup\{t_j;j=1...k\}$, and
note that $|\sum_{j=1}^k \theta_j N_{nt_j}(x)|\leq C N_{nT}(x)$ and
also $n[\tau n^H]^{-1}d(l,n)\leq C[\tau n^H]^{-1}\sum_{l\tau n^H\leq
  y< (l+1)\tau n^H}N_{nT}(y)$. Using Hölder and Minkowski's
inequalities, the second factor is  bounded from above by 
\begin{eqnarray*}
& & \left[\sum_{|l|\leq N}\sum_{l\tau n^H\leq x< (l+1)\tau n^H}\bbE \left( |\sum_{j=1}^k \theta_j N_{nt_j}(x)|^{\beta-1}+|n[\tau n^H]^{-1}d(l,n)|^{\beta-1}\right)^2 \right]^{1/2}\\
&\leq & C \left[\sum_{|l|\leq N}\sum_{l\tau n^H\leq x< (l+1)\tau n^H}\bbE \left(  N_{nT}(x)\right)^{2\beta-2} \right]^{1/2}\\
& & + C \left[\sum_{|l|\leq N}\sum_{l\tau n^H\leq x< (l+1)\tau
    n^H}\bbE \left(  [\tau n^H]^{-1}\sum_{l\tau n^H\leq y< (l+1)\tau
      n^H}N_{nT}(y) \right)^{2\beta-2}\right]^{1/2}
\end{eqnarray*}
\begin{eqnarray*}
&\leq &  C \left[\sum_{|l|\leq N}\sum_{l\tau n^H\leq x< (l+1)\tau n^H}\left(\bbE  N_{nT}(x)^2\right)^{\beta-1} \right]^{1/2}\\
& & + C \left[\sum_{|l|\leq N}\sum_{l\tau n^H\leq x< (l+1)\tau n^H}[\tau n^H]^{2-2\beta}\left(\sum_{l\tau n^H\leq y< (l+1)\tau n^H} \left(\bbE N_{nT}(y)^2   \right)^{1/2}\right)^{2\beta-2}\right]^{1/2}\\
&\leq & C\left[(2N+1)\tau n^H n^{(2-2H)(\beta-1)}\right]^{1/2}
\end{eqnarray*}
where the last line follows from Lemma 4.4 in \cite{W03} stating that there is some $C>0$ such that
 $$\sup_{x\in\bbsZ} \bbE(N_{nT}(x))^2\leq C n^{2-2H}.$$
Combining these estimates,  
 \begin{eqnarray*}
 & &\bbE \left| \sum_{|l|\leq N} \sum_{l\tau n^H\leq x< (l+1)\tau n^H}n^{-\delta\beta}\left\{ \big|\sum_{j=1}^k \theta_j N_{nt_j}(x)\big|^\beta -n^\beta[\tau n^H]^{-\beta}|d(l,n)|^\beta\right\} \right|\\
&\leq& Cn^{-\delta\beta}[(2N+1)\tau n^H\tau^rn^{2-2H}]^{1/2}\left[(2N+1)\tau n^H n^{(2-2H)(\beta-1)}\right]^{1/2}\\
&\leq & C N\tau^{1+r/2}
\end{eqnarray*}
Then the proof of Lemma \ref{KS} follows from equation (\ref{eqq}) and
from the above estimates as in \cite{KS79} or \cite{W03}: the idea is
to show that for large $N$, $n$ and small $\tau$, $X_n- V(\tau,n,N)
\rightarrow 0$ in probability and that $V(\tau,n,N)\rightarrow  X$ in
distribution. We omit the details. 

Next we  prove that $(X_n)_{n \ge 1} $ is bounded $L^p$ bound.
Let $T=\sup(t_1,\cdots,t_n)$ and $\Theta=\sum_{j=1}^k |\theta_j|$. 
The random variables $|X_n|$ is bounded above by 
$$\Theta^\beta n^{-\delta\beta} \sum_{x\in\bbsZ} N_{[nT]+1}^\beta(x).$$
In the case $\beta=2$, this quantity is equal to $\Theta^2 n^{H-2} V_{[nT]+1}$, and in this case the $L^p$ bound is a consequence of equation (\ref{eq3.5}).\\
In the case $\beta<2$,  Hölder inequality yields
$$\sum_{x\in\bbsZ} N_{[nT]+1}^\beta(x)\leq  \left(\sum_{x\in\bbsZ} 1_{\{N_{[nT]+1}(x)\neq0\}}\right) ^{1-\frac{\beta}{2}}\left(\sum_{x\in\bbsZ} N_{[nT]+1}^2(x)\right)^{\frac{\beta}{2}} =R_{[nT]+1}^{1-\frac{\beta}{2}}V_{[nT]+1}^{\frac{\beta}{2}}.$$
Hence, up to a multiplicative constant, the expectation $\bbE(|X_n|^p)$ is overestimated by
$$\bbE\left[ \left(n^{-\delta\beta} \sum_{x\in\bbsZ} N_{[nT]+1}^\beta(x)\right)^p\right]\leq 
\bbE \left[ \left(n^{-H}R_{[nT]+1} \right)^{p(1-\frac{\beta}{2})}\left(n^{H-2}V_{[nT]+1}\right)^{p\frac{\beta}{2}}\right]$$
We now apply Cauchy-Schwartz inequality,  
\begin{eqnarray*}
& &\bbE\left[ \left(n^{-\delta\beta} \sum_{x\in\bbsZ} N_{[nT]+1}^\beta(x)\right)^p\right]\\
&\leq& \bbE \left[ \left(n^{-H}R_{[nT]+1} \right)^{p(2-\beta)}\right]^{\frac{1}{2}} \bbE\left[\left(n^{H-2}V_{[nT]+1}\right)^{p\beta}\right]^{\frac{1}{2}}.
\end{eqnarray*}
Now the $L^p$ bound follows from equation (\ref{eq3.5}) and (\ref{eq3.6}) together.
$\Box$

\subsection{Convergence of the finite-dimensional distributions.}
We study the asymptotic behaviour of the characteristic function of the marginals of $G_n$. 
Let $\lambda$ be the characteristic function of the variables $\xi_k^{(i)}$ defined by
$$\lambda(u)=\bbE\left(\exp(iu\xi_1^{(1)})\right).$$
Since the random variables $\xi_k^{(i)}$ are in the domain of attraction of $Z_{\beta}$, 
\begin{equation}\label{eq3.15}
\lambda(u)=\bar\lambda (u)+o(|u|^\beta)\ \ ,\ \ {\rm as}\ u\rightarrow 0,
\end{equation}
where $\bar\lambda$ is the characteristic function of $Z_\beta$ given by equation (\ref{eq1.02}).

\begin{pr}\label{pr2.1}
The finite dimensional distributions of $(G_n(t))_{t\geq 0}$ converge weakly as $n\rightarrow\infty$ to those of $(\Gamma(t))_{t\geq 0}$ defined in equation (\ref{eq1.09}).
\end{pr}

{\it Proof:}\\
Let $(\theta_1,\cdots,\theta_k)\in \bbR^k$ , $(t_1,\cdots,t_k)\in [0,+\infty)^k$. Computations as in \cite{DG07} show that the characteristic function of $\Gamma_n(t)$ writes
\begin{equation}\label{eq4.1}
\bbE\left[\exp\left(i\sum_{j=1}^k \theta_j G_n(t_j)\right) \right]
= \left(\bbE \left[ \prod_{x\in\bbsZ} \lambda \left( c_n^{-\frac{1}{\beta}}U_n(x)\right)\right]\right )^{c_n} 
\end{equation}
where
$$U_{n}(x)=n^{-\delta} \sum_{j=1}^k \theta_j N_{nt_j}(x)\ \ ,\ \ x\in\bbZ.$$
We show that the following asymptotic holds as $n\rightarrow \infty$:
\begin{equation}\label{eq4.2}
  \bbE\left[ \prod_{x\in\bbsZ} \lambda \left( c_n^{-\frac{1}{\beta}}U_{n}(x)\right) \right] = \bbE\left[ \prod_{x\in\bbsZ} \lambda \left( c_n^{-\frac{1}{\beta}}U_{n}(x)\right) \right] +o(c_n^{-1})
\end{equation}
To see this, note that
\begin{eqnarray}
 & &c_n\left|\ \prod_{x\in\bbsZ} \lambda \left( c_n^{-\frac{1}{\beta}}U_{n}(x)\right)- \prod_{x\in\bbsZ} \bar \lambda \left( c_n^{-\frac{1}{\beta}}U_{n}(x)\right)\ \right| \nonumber \\
 &\leq& c_n \sum_{x\in\bbsZ}  \ \left|\ \lambda \left( c_n^{-\frac{1}{\beta}}U_{n}(x)\right)- \bar\lambda \left( c_n^{-\frac{1}{\beta}}U_{n}(x)\right)\  \right|\nonumber \\
&\leq & \tilde g(c_n^{-\frac{1}{\beta}}U_{n}) \sum_{x\in\bbsZ}|U_{n}(x)|^\beta \label{eq3.17}.
\end{eqnarray}
with 
$$U_{n}=\sup_{x\in\bbsZ} |U_{n}(x)|,$$
and $\tilde g$ the bounded continuous vanishing at zero function defined by
$$\tilde g(u)=\sup_{|v|\leq u}
|v|^{-\beta}\left|\lambda(v)-\bar\lambda(v)\right|\ \ ,\ \ v\neq 0.$$
(The properties of $\tilde g$ follow from equation (\ref{eq3.15}).)
From Lemma \ref{lem3.1}, $U_{n}$ converge in probability to $0$ as $n\rightarrow\infty$. Since $\tilde g$ is bounded continuous and vanishes at $0$, $\tilde g(c_n^{-\frac{1}{\beta}}U_{n})$ converges also in probability to $0$ and is bounded in $L^\infty$. 
From Lemma \ref{KS}, $\sum_{x\in\bbsZ}|U_{n}(x)|^\beta$ converges in distribution and is bounded in $L^p$.  As a consequence, the right hand side of (\ref{eq3.17}) converges to zero in probability and is bounded in $L^p$, and hence its expectation has limit $0$. This proves equation (\ref{eq4.2}).\\

We now prove the following estimation
\begin{equation}\label{eq4.3}
  \bbE\left[ \prod_{x\in\bbsZ} \bar \lambda \left( c_n^{-\frac{1}{\beta}}U_{n}(x)\right) \right] = 1- c_n^{-1}\sigma^\beta\bbE\left[ X \right] +o(c_n^{-1})
\end{equation}
where $X$ is defined in Lemma \ref{KS}.To see this, recall the definition of the random variable $X_n$ from Lemma \ref{KS} and of the characteristic function $\bar\lambda$ from equation (\ref{eq1.02}). With these notations, equation (\ref{eq4.3}) is equivalent to 
$$\lim_{n\rightarrow +\infty}  \bbE\left(f_n(X_n) \right)=\sigma^\beta \bbE(X),$$
where $f_n$ is the function defined on $\bbC$ by 
$$f_n(x)=c_n\left(1-\exp(-c_n^{-1}\sigma^\beta x)  \right).$$
It is easy to verify that the sequence of functions $f_n$ satisfies the following property: for every $x$, for every sequence $(x_n)_{n\geq 1}$ converging to $x$, 
$$\lim_{n\rightarrow\infty} f_n(x_n)=\sigma^\beta x.$$
Furthermore Lemma \ref{KS} states that the sequence $(X_n)_{n\geq 1}$ converges in distribution to $X$ when $n\rightarrow\infty$. Using the diagonal mapping Theorem (Theorem 5.5 of \cite{Bi68}), we prove the weak convergence of the sequence of random variables $f_n(X_n)$ to $\sigma^\beta X$. Furthermore, using Lemma \ref{KS}, $\left|f_n(X_n)\right|\leq |X_n|$ is bounded in $L^p$ for any $p\geq 1$. Hence $\bbE(f_n(X_n))$ has limit $\sigma^\beta \bbE(X)$ and equation (\ref{eq4.3}) is proved.\\
Finally, combining equations (\ref{eq4.1}), (\ref{eq4.2}) and (\ref{eq4.3}) we prove easily that
$$ \bbE\left[\exp\left(i\sum_{j=1}^k \theta_j G_n(t_j)\right) \right] =\left(1- c_n^{-1}\ \sigma^\beta\bbE(X) +o(c_n^{-1})\right)^{c_n}\mathop{\longrightarrow}_{n\rightarrow\infty}\exp(-\sigma^\beta\bbE(X))$$
and Proposition \ref{pr2.1} is proved.
$\Box$

\subsection{Tightness}
\begin{pr}\label{pr2.2}
The family of processes $\left(G_{n}(t)\right)_{t\geq 0}$ is tight in $\cC([0,\infty))$.
\end{pr}
{\it Proof:}
As in the continuous case, we prove the tightness using truncations in order to deal with finite variance processes.  We decompose the scenery $(\xi_x^{(i)})_{x\in\bbsZ,i\geq 1}$ into two parts
$$\xi_x^{(i)}=\bar\xi_{a,x}^{(i)}+\hat\xi_{a,x}^{(i)},$$
where $(\bar\xi_{a,x}^{(i)})$ denote the $i$-th truncated scenery defined by
$$\bar\xi_{a,x}^{(i)}=\xi_x^{(i)} 1_{\{ |\xi_x^{(i)}|\leq a\}},$$
and $\hat\xi_{a,x}^{(i)}$ the remainder scenery  
$$\hat\xi_{a,x}^{(i)}=\xi_x^{(i)} 1_{\{ |\xi_x^{(i)}| > a\}}.$$
We recall the following estimates from Lemma 3.3 in \cite{DG07}: there exists some $C>0$ such that
\begin{equation}\label{eq3.22}
 |\bbE\left(\bar\xi_{a,x}^{(i)}\right)|\leq Ca^{1-\beta}\ ,\ \bbE\left(|\bar\xi_{a,x}^{(i)}|^2\right)\leq  Ca^{2-\beta}\ ,\ 
\bbP\left(\hat \xi_{a,x}^{(i)} \neq 0 \right) \leq   C a^{-\beta}. 
\end{equation}
For $a>0$, we use truncations with $a_n=a n^{\frac{H}{\beta}}c_n^{\frac{1}{\beta}}$ and write
\begin{equation}\label{eq3.20}
G_n(t)=\bar \Gamma_{n,a}(t)+\hat\Gamma_{n,a}(t),
\end{equation}
where
\begin{eqnarray*}
 \bar \Gamma_{n,a}(t)&=& n^{-\delta}c_n^{-\frac{1}{\beta}}\sum_{i=1}^{c_n}  \sum_{x\in\bbsZ} N_{nt}^{(i)}(x) \bar\xi_{a_n,x}^{(i)},\\
 \hat \Gamma_{n,a}(t)&=& n^{-\delta}c_n^{-\frac{1}{\beta}}\sum_{i=1}^{c_n}  \sum_{x\in\bbsZ} N_{nt}^{(i)}(x) \hat\xi_{a_n,x}^{(i)}.
\end{eqnarray*}
Now, with the same techniques as in the proof of Proposition 3.2 in \cite{DG07}, we compute:
\begin{eqnarray*}
\bbP\left( \sup_{t\in[0,T]}|\hat \Gamma_{n,a}(t)|=0\right) &\geq & \left(\bbE\left[\left(\bbP (\hat\xi_{a_n,0}^{(1)}=0)\right)^{R_{[nT]+1}}\right]\right)^{c_n}\\
 &\geq& \left(\bbE \left[\left(1-Ca_n^{-\beta} \right)^{R_{[nT]+1}}\right]\right)^{c_n}\\
 &\geq & \left(1+ \log(1-Ca_n^{-\beta}) \bbE( R_{[nT]+1})\right)^{c_n}
\end{eqnarray*}
Using the asymptotic for $a_n$ and Lemma \ref{lem3.1} to estimate the range, the above inequality implies 
\begin{equation}\label{eq3.19}
\lim_{a\rightarrow \infty}\limsup_{n\rightarrow\infty}\bbP\left( \sup_{t\in[0,T]}|\hat \Gamma_{n,a}(t)|>0\right)=0
\end{equation}
On the other hand, the variance of the truncated process $\Gamma_{n,a}$ is overestimated by
\begin{eqnarray*}
& &\bbE\left[\left|\ \bar \Gamma_{n,a}(t_2)-\bar \Gamma_{n,a}(t_1)\ \right|^2\right]\nonumber\\
&\leq&n^{-2\delta}c_n^{-\frac{2}{\beta}}c_n(c_n-1) \left[\sum_{x\in\bbsZ} \bbE(N_{nt_2}^{(1)}(x)-N_{nt_1}^{(1)}(x))\right]^2\left[\bbE|\bar\xi_{a_n,0}^{(1)}|  \right]^2\\
& &+n^{-2\delta}c_n^{-\frac{2}{\beta}}c_n \sum_{x\neq y\in\bbsZ}\bbE\left[(N_{nt_2}^{(1)}(x)-N_{nt_1}^{(1)}(x))(N_{nt_2}^{(1)}(y)-N_{nt_1}^{(1)}(y))\right]\left[\bbE|\bar\xi_{a_n,0}^{(1)}| |\right]^2\\
& &+n^{-2\delta}c_n^{-\frac{2}{\beta}}c_n \sum_{x\in\bbsZ}\bbE\left[(N_{nt_2}^{(1)}(x)-N_{nt_1}^{(1)}(x))^2\right]\bbE\left[|\bar\xi_{a_n,0}^{(1)}|^2 \right]\\
\end{eqnarray*}
Using equation (\ref{eq3.22}) and the following estimations, 
$$\bbE\left[ \sum_{x\in\bbsZ} (N_{nt_2}^{(1)}(x)-N_{nt_1}^{(1)}(x))\right]=n(t_2-t_1) $$
$$\bbE\left[ \sum_{x\neq y\in\bbsZ} (N_{nt_2}^{(1)}(x)-N_{nt_1}^{(1)}(x)) (N_{nt_2}^{(1)}(y)-N_{nt_1}^{(1)}(y))\right]=n^2(t_2-t_1)^2- \bbE\left[\sum_{x\in\bbsZ}(N_{nt_2}^{(1)}(x)-N_{nt_1}^{(1)}(x))^2\right]$$
$$ \bbE\left[\sum_{x\in\bbsZ} (N_{nt_2}^{(i)}(x)-N_{nt_1}^{(i)}(x))^2 \right] \leq \bbE(V_{[nt_2]-[nt_1]+1}) \leq C([nt_2]-[nt_1]+1)^{2-H}$$
we prove that there exists some $C$ such that if $|t_2-t_1|\geq \frac{1}{n}$, then
$$\bbE\left[\left|\ \bar \Gamma_{n,a}(t_2)-\bar \Gamma_{n,a}(t_1)\ \right|^2\right]\leq C|t_2-t_1|^{2-H}.$$
In the case   $|t_2-t_1|\leq 1/n$, we can see that  
$$\bbE\left[\sum_{x\in\bbsZ} (N_{nt_2}^{(1)}(x)-N_{nt_1}^{(1)}(x))^2 \right]\leq 2 (nt_2-nt_1)^2,$$
since in the sum, at most 
two  terms are not zero and those terms are bounded by $(nt_2-nt_1)^2$.
Using theorem 12.3 in Billingsley, these estimates prove the tightness of the family of processes $\left(\bar \Gamma_{n,a}(t)\right)_{t\geq 0}$. This together with equations (\ref{eq3.19}) and (\ref{eq3.20}) implies the tightness of the sequence $G_n$, and hence Proposition \ref{pr2.2}.

$\Box$

%


\end{document}